# A GENERATING FUNCTION FOR THE TRACE OF THE IWAHORI-HECKE ALGEBRA

ERIC M. OPDAM

ABSTRACT. The Iwahori-Hecke algebra has a "natural" trace $\tau$. This trace is the evaluation at the identity element in the usual interpretation of the Iwahori-Hecke algebra as a sub-algebra of the convolution algebra of a p-adic semi-simple group. The Iwahori-Hecke algebra contains an important commutative sub-algebra $\mathbf{C}[\theta_x]$, that was described and studied by Bernstein, Zelevinski and Lusztig. In this note we compute the generating function for the value of $\tau$ on the basis $\theta_x$.

*Date*: October 25, 2018.

1991 *Mathematics Subject Classification.* 20C08, 22D25, 22E35.

The author would like to thank Erik van den Ban, Gerrit Heckman and Klaas Slooten for many useful remarks and fruitful discussions.

CONTENTS

## 1. Notations and main result

The goal of this paper is to derive a disintegration formula which could serve as a starting point for the spectral analysis of the Iwahori-Hecke algebra. However, I have decided to postpone the treatment of the more spectral analytic aspects of this formula to another paper. The formula has a certain interest in its own right, and as a first application we will give an explicit formula for the value of the trace on the maximal Abelian subalgebra of the Hecke-algebra.

I hope to help the reader (and myself) by reviewing the setup and definition of the Hecke-algebra, and the elementary theory of its minimal principal series representations. This is mainly to fix the notations, but also to give the reader an impression of the logical order of the basic facts that are involved in the proof of the main formula of this paper. At the same time, this hopefully provides a convenient basis for further papers in this direction. As a drawback, the material in the first section, and parts of the second section is (very) classical and well known to specialists.

In the first section we review the definition of the Iwahori-Hecke algebra. We basically follow Lusztig's paper [3] in this discussion.

## 1.1. The root datum

Let $\mathcal{R} = (X, Y, R_0, R_0^\vee, F_0)$ be a reduced root datum (I apologize for the subscripts, but I have decided to reserve the notation $R$ for the affine root system which has to be introduced later). This means that

(a) $X$ and $Y$ are lattices, with a given perfect pairing
$$(\cdot, \cdot) : X \times Y \to \mathbf{Z},$$

(b) $R_0 \subset X$, $R_0^\vee \subset Y$ are finite subsets with a given bijection $R_0 \to R_0^\vee$, denoted by $\alpha \to \alpha^\vee$. (The elements of $R_0$ are called "roots", those of $R_0^\vee$ are called "coroots".)

(c) $F_0 \subset R_0$ (called the set of fundamental roots).

(d) $(\alpha, \alpha^\vee) = 2$ for all $\alpha \in R_0$.

(e) For any $\alpha \in R_0$, the reflection $s_\alpha \in GL(X)$ given by $s_\alpha(x) = x - (x, \alpha^\vee)\alpha$ maps $R_0$ to itself. Similarly, $s_\alpha \in GL(Y)$ given by $s_\alpha(y) = y - (\alpha, y)\alpha^\vee$ maps $R_0^\vee$ to itself.

(f) $F_0$ is linearly independent, and $R_0 \subset \mathbf{Z}_{\geq 0} F_0 \cup \mathbf{Z}_{\leq 0} F_0$. (The set $R_{0,+} = R_0 \cap \mathbf{Z}_{\geq 0} F_0$ is called the set of positive roots.)

(g) For all $\alpha \in R_0$, $2\alpha \notin R_0$.

Put $Q = \mathbf{Z} R_0$ and $Q^\vee = \mathbf{Z} R_0^\vee$ for the root lattice and the coroot lattice respectively. The weight lattice $P$ of $R$ is by definition the dual $P = \text{Hom}(Q^\vee, \mathbf{Z})$. It follows from the definitions that there is a natural map $X \to P$ whose restriction to $Q \subset X$ is the natural inclusion $Q \subset P$.

The group generated by the reflections $s_\alpha$ is called the Weyl group of $\mathcal{R}$, denoted $W_0$. The set $S_0 = \{s_\alpha \mid \alpha \in F\}$ is called the set of simple reflections. It is well known that $S_0$ generates $W_0$, and that $(W_0, S_0)$ is a finite Coxeter group.

We define the so called dominance (partial) ordering on $Y$ by:
$$y \geq y' \Leftrightarrow y - y' \in \mathbf{Z}_{\geq 0} R_{0,+}^\vee.$$

## 1.2. The affine Weyl group

The affine Weyl group $W$ is the semi-direct product
$$W = W_0 \ltimes X.$$

The elements of $W$ are denoted $wt_x$, with the multiplication rule $wt_x \cdot w't_{x'} = ww't_{w'^{-1}x + x'}$. Clearly $W$ acts on the set $X$. When we think of $Y \times \mathbf{Z}$ as a set of affine linear functions on $X$, we get the natural dual action of $W$ on $Y \times \mathbf{Z}$, defined by $(wf)(x) := f(w^{-1}x)$. Explicitly, this amounts to:
$$wt_x(y, k) := (w(y), k - (x, y)).$$
The affine root system $R$ is by definition the set $R_0^\vee \times \mathbf{Z} \subset Y \times \mathbf{Z}$. It is the a $W$-invariant set in $Y \times \mathbf{Z}$ containing the set of coroots $R_0^\vee$. It is easy to see that $R = R_+ \cup R_- = R_+ \cup -R_+$ where:
$$R_+ = \{(\alpha^\vee, k) \mid k > 0 \text{ and } \alpha \in R_0\} \cup \{(\alpha^\vee, 0) \mid \alpha^\vee \in R_{0,+}^\vee\}.$$
Let $S^m \subset R_0$ be the set of roots so that $\{S^m\}^\vee$ is the set of minimal coroots (with respect to the dominance ordering). It is clear that $R_+$ is a simplicial cone, generated over $\mathbf{Z}_{\geq 0}$ by the set
$$F = \{(\alpha^\vee, 1) \mid \alpha \in S^m\} \cup \{(\alpha^\vee, 0) \mid \alpha \in F_0\}.$$
This is called the set of fundamental affine roots. For every affine root $a = (\alpha^\vee, k)$ we define a corresponding affine reflection $s_a : X \to X$ by the formula
$$s_a(x) = x - a(x)\alpha,$$
where $a(x) = (x, \alpha^\vee) + k$. This is an element of $W$ since $s_a = s_\alpha t_{k\alpha}$. We put
$$S = \{s_a \mid a \in F\},$$
the set of affine simple reflections.

It follows easily that the subgroup of $W$ generated by $S$ equals $W_0Q$. This is a normal subgroup in $W$, which is an affine Coxeter group with $S$ as set of simple reflections. We shall sometimes think of $S$ as the set of vertices of the affine Coxeter diagram of $(W_0Q, S)$.

### 1.3. The length function

The length function $l$, defined on elements $w$ of $W_0Q$ as the minimal length of a word in the alphabet $S$ representing $w$, has a natural extension to $W$. The extension is defined as follows:
(1.1)
$$\begin{aligned} l(wt_x) &= \#\{a \in R_+ \mid wt_x(a) \in R_-\} \\ &= \sum_{\alpha \in R_{0,+} \cap w^{-1}(R_{0,-})} |(x, \alpha^\vee) + 1| + \sum_{\alpha \in R_{0,+} \cap w^{-1}(R_{0,+})} |(x, \alpha^\vee)|. \end{aligned}$$
Indeed, the length of an element $s_a \in S$ corresponding to $a \in F$ is equal to 1, since $R_+ \cap s_a(R_-) = \{a\}$. It follows that for all $w \in W$,
$$(1.2) \qquad l(s_a w) = \begin{cases} l(w) + 1 & \text{if } w^{-1}(a) \in R_+. \\ l(w) - 1 & \text{if } w^{-1}(a) \in R_-. \end{cases}$$

From 1.2 it is standard to derive that the restriction of $l$ to $W_0 Q$ indeed coincides with the length function of this Coxeter group with respect to the set of generators $S$.

Put $X^+ = \{x \in X \mid (x, \alpha^\vee) \geq 0, \forall \alpha \in F_0\}$, the set of dominant elements of $X$. Observe that 1.1 implies that for all $x \in X^+$, $w \in W_0$:

$$\begin{aligned} l(wt_x) &= l(w) + l(t_x) \\ &= l(w) + 2(x, \rho^\vee) \end{aligned} \tag{1.3}$$

where $\rho^\vee \in P^\vee$ is defined by

$$\rho^\vee = \frac{1}{2} \sum_{\alpha \in R_{0,+}} \alpha^\vee.$$

Similarly, when $x \in X^{++}$ is strongly dominant and $w \in W_0$ we have the formula

$$l(t_x w) = -l(w) + l(t_x).$$

**Corollary 1.1.** *The length is an additive function on the Abelian half-group consisting of $t_x$, $x \in X^+$.*

*Proof.* Clear from 1.3. □

It is easy to see from the above description of the length function that an element $w \in W$ has length 0 if and only if it maps $F$ to itself. Thus the set of elements in $W$ having length 0 is a subgroup $\Omega$ of $W$. Using 1.2 it is clear that this is a complement for the normal Coxeter subgroup $W_0 Q$. Hence $\Omega$ is Abelian, and isomorphic to the quotient $X/Q$.

### 1.4. Root labels

In order to define the Iwahori-Hecke algebra we need to fix a length-multiplicative function $q$ on $W$. In principle $q$ may take its values in any Abelian half-group, but for our purposes it is appropriate to assume that $q$ takes its values in $\mathbf{R}_{>1}$. In other words, we consider a function

$$q : W \to \mathbf{R}_{>1}$$

such that $q(ww') = q(w)q(w')$ whenever $l(ww') = l(w) + l(w')$.

Note that $q$ is identically equal to 1 on the group $\Omega$ of elements of length 0. The function $q$ is therefore completely determined by its values on the set $S$ of elements of length 1. On the other hand, any function $q$ on $S$ can be extended length-multiplicatively provided that $q_s = q_{s'}$ whenever $s$, $s' \in S$ are conjugate in $W$.

We associate to such a function $q$ on $W$ a $W$-invariant set of root labels on the affine root system $R$ by the following rules. First observe that a translation $t_x$ with $x \in X$ maps $a = \alpha^\vee + k \in R$ to $a - (x, \alpha^\vee)$. Hence $a$ is in the orbit of $\alpha^\vee$, except when $\alpha^\vee \in 2Y$ and $k$ is odd. Let $a = \alpha^\vee + k \in F$. If $\alpha^\vee \notin 2Y$, we simply put

$$q_a := q(s_\alpha).$$

When $\alpha^\vee \in 2Y$, let $S(a) \subset S$ denote the connected component of the Coxeter graph of the affine Coxeter group $(W_0Q, S)$. As Lusztig ([3], Lemma 1.7) remarks, $S(a)$ is of affine type $C_n^{\text{aff}}$, and the corresponding root datum $\mathcal{R}(a)$ of type $(Q(B_n) = \mathbf{Z}^n, P(C_n) = \mathbf{Z}^n, B_n, C_n, \{e_1 - e_2, \ldots, e_{n-1} - e_n, e_n\})$ is a direct summand of $\mathcal{R}$. Let $\omega_a$ be the unique nontrivial diagram automorphism of $S$ which is trivial on the complement of $S(a)$. We now define

$$q_a := q(s_{\omega_a(a)}).$$

Explicitly, we define for these fundamental roots in $S(a)$: $q_{1-2e_1} = q(s_{2e_n})$ and $q_{2e_n} = q(s_{1-2e_1})$. We extend the affine root labels $W$-invariantly to all of $R$. It is clear that the affine root labels determine $q$.

Finally, we can also describe $q$ by means of a $W_0$-invariant set of root labels $q_{\alpha^\vee}$, where $\alpha$ runs over the roots of the possibly non-reduced root system $R_{\text{nr}}$ defined by

$$R_{\text{nr}} := R_0 \cup \{2\alpha \mid \alpha^\vee \in R_0^\vee \cap 2Y\}.$$

We do this by restricting the affine root labels to $R_0^\vee$, and when $\alpha^\vee \in 2Y$ we define:

$$q_{\alpha^\vee/2} := \frac{q_{1+\alpha^\vee}}{q_{\alpha^\vee}}.$$

It is clear that this set of root labels on $R_{\text{nr}}^\vee$ also contains the same information as the function $q$. These definitions of root labels are in compliance with Macdonald [4], except for the fact that we call roots what he calls coroots!

We list some direct consequences of our notations:

**Corollary 1.2.** *For all $w \in W$,*

$$q(w) = \prod_{a \in R_+ \cap w^{-1} R_-} q_{a+1}.$$

**Corollary 1.3.** *For all $w \in W_0$,*

$$q(w) = \prod_{\alpha \in R_{\text{nr},+} \cap w^{-1} R_{\text{nr},-}} q_{\alpha^\vee}.$$

**Definition 1.4.** *The Haar modulus $\delta = \delta_q$ is the character on on $X$ given by*
$$\delta(x) = \delta_q(x) := \prod_{\alpha \in R_{\mathrm{nr},+}} q_{\alpha^\vee}^{(x,\alpha^\vee)}.$$

**Corollary 1.5.** *When $x \in X^+$,*
$$q(t_x) = \delta(x).$$

### 1.5. The Iwahori-Hecke algebra

Now we are ready to define the Iwahori-Hecke algebra $\mathcal{H}$ associated to a root datum $\mathcal{R}$ and a set of root labels $q_{\alpha^\vee}$ with $\alpha \in R_{\mathrm{nr}}$.

**Definition 1.6.** *The Iwahori-Hecke algebra $\mathcal{H} = \mathcal{H}(\mathcal{R}, q)$ is the complex associative algebra with basis $T_w$, $w \in W$, and relations induced by the following rules:*
  (a) *If $l(ww') = l(w) + l(w')$ then $T_w T_{w'} = T_{ww'}$.*
  (b) *If $s \in S$ then $(T_s + 1)(T_s - q(s)) = 0$.*

**Remark 1.7.** *It is well known that such an object exists.*

This algebra can be equipped with a $*$ operator, which is the anti-linear anti-involution defined by
$$T_w^* = T_{w^{-1}}.$$
It is straightforward to check that this indeed extends to an anti-involution on $\mathcal{H}$.

The object of study in this note is the canonical trace function $\tau$, which is by definition the linear functional on $\mathcal{H}$ that has the following values on the basis elements $T_w$:
$$\tau(T_w) = \begin{cases} 1 & \text{if } w = e. \\ 0 & \text{else.} \end{cases}$$

This functional corresponds to the evaluation at the identity element when the Hecke algebra is interpreted as a subalgebra of the convolution algebra of a semi-simple group of $p$-adic type. Therefore the following lemma is natural:

**Proposition 1.8.** *The trace $\tau$ is central and positive.*

*Proof.* Both statements are clear from the following obvious formula:
(1.4) $$\tau(T_w^* T_{w'}) = \delta_{w,w'} q(w).$$

□

The associated Hermitian inner product on $\mathcal{H}$ is given by the formula
$$(h, h') := \tau(h^* h').$$
Observe that the basis $T_w$ is orthogonal with respect to this inner product because of 1.4.

## 1.6. The Bernstein-Zelevinski basis

The Hecke algebra $\mathcal{H}$ contains a large commutative subalgebra. For $x \in X^+$ we define
$$\theta_x = \delta(-x)^{1/2} T_{t_x}.$$
These elements form an Abelian half-group isomorphic to $X^+$ because of Corollary 1.1. By 1.5 they are orthonormal. For general $x \in X$ we define:
$$\theta_x = \theta_y \theta_z^{-1}$$
where $y, z \in X^+$ are such that $x = y - z$.

The following result is very important, and due to Bernstein and Zelevinski (unpublished work).

**Theorem 1.9.**  (a) *The elements $T_w \theta_x$, $w \in W_0$, $x \in X$ form a basis of $\mathcal{H}$.*
  (b) *The elements $\theta_x T_w$, $w \in W_0$, $x \in X$ also form a basis of $\mathcal{H}$.*
  (c) *In particular, the subalgebra $\mathcal{A} = \mathbf{C}[\theta_x] \subset \mathcal{H}$ is isomorphic to the group algebra of $X$ via $x \to \theta_x$.*

The next result is due to Bernstein and Zelevinski in a special case, and to Lusztig in the general case. It tells us how the multiplication works in terms of the basis just described in the above theorem.

**Theorem 1.10.** *Let $x \in X$ and $\alpha \in F_0$. Let $s = s_\alpha$. Then*

(1.5)
$$\theta_x T_s - T_s \theta_{s(x)} =$$
$$= \begin{cases} (q_{\alpha^\vee} - 1) \frac{\theta_x - \theta_{s(x)}}{1 - \theta_{-\alpha}} & \text{if } 2\alpha \notin R_{\text{nr}}. \\ ((q_{\alpha^\vee/2} q_{\alpha^\vee} - 1) + q_{\alpha^\vee/2}^{1/2}(q_{\alpha^\vee} - 1)\theta_{-\alpha}) \frac{\theta_x - \theta_{s(x)}}{1 - \theta_{-2\alpha}} & \text{if } 2\alpha \in R_{\text{nr}}. \end{cases}$$

One important consequence of this theorem is the precise description of the center $\mathcal{Z}$ of $\mathcal{H}$. This result is also due to Bernstein and Zelevinski.

**Theorem 1.11.** *The subalgebra $\mathcal{A} \subset \mathcal{H}$ is a $W_0$ module via the action of $W_0$ on $X$. In this sense, the center $\mathcal{Z}$ equals $\mathcal{Z} = \mathcal{A}^{W_0}$.*

The behavior of the elements $\theta_x$ with respect to the $*$ operator is not very complicated:

**Proposition 1.12.** *Let $w_0 \in W_0$ denote the longest element.*
$$\theta_x^* = T_{w_0}\theta_{-w_0(x)}T_{w_0}^{-1}.$$

*Proof.* We may assume that $x \in X^+$. We need to show that $T_{t_{-x}}T_{w_0} = T_{w_0}T_{t_{-w_0(x)}}$. This is obvious since $l(w_0 t_{-w_0(x)}) = l(w_0) + l(t_x)$ (see 1.3) and $w_0 t_{-w_0(x)} = t_{-x}w_0$. □

### 1.7. The main result

Let $T = \text{Hom}(X, \mathbf{C}^\times)$ be the complex torus of characters of $X$. Let us introduce the $c$-function, following Macdonald [4]. For $\alpha \in R_{\text{nr}}$ and $t \in T$ define
$$c(\alpha, t) := \frac{1 - q_{\alpha^\vee/2}^{-1/2}q_{\alpha^\vee}^{-1}t(-\alpha)}{1 - q_{\alpha^\vee/2}^{-1/2}t(-\alpha)}.$$

(In this notation it is understood that $q_\beta^\vee = 1$ whenever $\beta \notin R_{\text{nr}}$, and that we take the positive root of $q_{\alpha^\vee/2}$). When $\alpha \in R_0$ we put
$$c_0(\alpha, t) = c(\alpha, t)c(2\alpha, t).$$

It is convenient to introduce the reduced root system $R_1 = \{\alpha \in R_{\text{nr}} \mid 2\alpha \notin R_{\text{nr}}\}$. We define for $\alpha \in R_1$:
$$c_1(\alpha, t) = c(\alpha, t)c(\alpha/2, t).$$

Clearly, when $\alpha \in R_0 \cap R_1$ we have $c_0(\alpha, t) = c_1(\alpha, t)$. On the other hand, when $\alpha \in R_0 - R_1$, then $2\alpha \in R_1$, and $c_0(\alpha, t) = c_1(2\alpha, t)$. We have for $\alpha \in R_0$:
$$c_0(\alpha, t) = \begin{cases} \frac{1 - q_{\alpha^\vee}^{-1}t(-\alpha)}{1 - t(-\alpha)} & \text{if } 2\alpha \notin R_{\text{nr}}. \\ \frac{(1 + q_{\alpha^\vee/2}^{-1/2}t(-\alpha))(1 - q_{\alpha^\vee/2}^{-1/2}q_{\alpha^\vee}^{-1}t(-\alpha))}{1 - t(-2\alpha)} & \text{if } 2\alpha \in R_{\text{nr}}. \end{cases}$$

With these notations we define:

**Definition 1.13.**
$$c(t) := \prod_{\alpha \in R_{\text{nr},+}} c(\alpha, t) = \prod_{\alpha \in R_{0,+}} c_0(\alpha, t) = \prod_{\beta \in R_{1,+}} c_1(\beta, t).$$

We introduce a partial ordering $<$ on the real split part $T_{rs} = \text{Hom}(X, \mathbf{R}_+)$ of $T$ by
$$t_1 < t_2 \Leftrightarrow \forall \alpha \in R_{0,+} : t_1(\alpha) < t_2(\alpha).$$

Here is the main result of this paper, a generating function for the values of $\tau$ on $\mathcal{A}$:

**Theorem 1.14.** *We have the following identity of formal Laurent series ($t \in T$):*

$$\sum_{x \in X} \tau(\theta_x) t(-x) = \bigl(q(w_0) c(t) c(t^{-1})\bigr)^{-1}_{\exp,+}.$$

*The right hand side of the equality means the expansion of the rational function $1/(q(w_0)c(t)c(t^{-1}))$ as a power series in $t_i = t(\alpha_i)$, where $F_0 = \{\alpha_1, \ldots, \alpha_n\}$. Notice that this expansion is convergent on $\{t \in T \mid \mathrm{Re}(t) < \delta^{-1/2}\}$.*

**Remark 1.15.** *In particular, this implies that $\tau(\theta_x) = 0$ unless $x \in Q_-$.*

**Example 1.16.** *The simplest possible case is when the root datum is of type $A_1$, i.e. $\mathcal{R} = (\mathbf{Z}, \mathbf{Z}, \{\pm 2\}, \{\pm 1\}, \{2\})$. Put $S_0 = \{s\}$, $q = q(s)$, and $\omega = t_1 s$. Notice that $\omega \in \Omega$, the group of elements of length 0 in $W$. Now consider $\theta = q^{-1/2} T_{t_1} = q^{-1/2} T_\omega T_s$. The above theorem states that $\tau(\theta^n) = 0$ unless $n = -2k$ with $k \in \mathbf{Z}_{\geq 0}$, and when $k > 0$ we have*

$$\tau(\theta^{-2k}) = \frac{(q-1)(q^k - q^{-k})}{q+1} = (q + q^{-1} - 2)(q^{k-1} + q^{k-3} + \cdots + q^{1-k}).$$

*The reader is invited to verify this formula directly.*

**Lemma 1.17.** *Let $\alpha \in R_{\mathrm{nr}}$ and $k \in \mathbf{Z}_{\geq 0}$. We define a rational function $d(\alpha; k)$ in $q^{1/2}_{\alpha^\vee/2}$ and $q_{\alpha^\vee}$ by*

$$\sum_{k=0}^{\infty} d(\alpha; k) t^k = \bigl(q_{\alpha^\vee} c(\alpha, t) c(\alpha, t^{-1})\bigr)^{-1}_{\exp,+}.$$

*Then $d(\alpha; 0) = 1$, and for $k > 0$ we have:*

$$d(\alpha; k) = \frac{(q_{\alpha^\vee} - 1)(q_{\alpha^\vee/2} q_{\alpha^\vee} - 1)\Bigl((q^{1/2}_{\alpha^\vee/2} q_{\alpha^\vee})^k - (q^{1/2}_{\alpha^\vee/2} q_{\alpha^\vee})^{-k}\Bigr)}{(q_{\alpha^\vee/2} q^2_{\alpha^\vee} - 1)}$$

**Corollary 1.18.** *We can reformulate the main result as follows. Let $x \in X$. We call $\pi = (\pi_\alpha)_{\alpha \in R_{\mathrm{nr},+}}$ a partition of $x$ if (1) $\pi_\alpha \in \mathbf{Z}_{\geq 0}$ and (2) $x = \sum_{\alpha \in R_{\mathrm{nr},+}} \pi_\alpha \alpha$. Then the trace of $\theta_x$ can be written as a weighted partition function:*

(1.6) $$\tau(\theta_x) = \sum_{\pi} \prod_{\alpha \in R_{\mathrm{nr},+}} d(\alpha; \pi_\alpha),$$

*where the sum is taken over all the partitions $\pi$ of $-x$.*

**Corollary 1.19.** *Suppose that $q(s) > 1$, $\forall s \in S$. Then $\tau(\theta_{-\kappa}) > 0$, for all $\kappa \in Q_+$.*

## 2. The minimal principal series and their Intertwining operators

### 2.1. Definition of the minimal principal series

The holomorphic minimal principal series are simply defined by induction from $\mathcal{A}$. Let $t \in T$ denote a character of $X$.

**Definition 2.1.** *The (holomorphic) minimal principal series $I_t$ is the induced module $I_t = Ind_{\mathcal{A}}^{\mathcal{H}}(t) = \mathcal{H} \otimes_{\mathcal{A}} \mathbf{C}_t$ (where $\mathbf{C}_t$ denotes the one dimensional module of $\mathcal{A}$ corresponding the character $t$).*

As an $\mathcal{H}_0$ module, $I_t \simeq \mathcal{H}_0$ via $T_w \otimes 1 \to T_w$ ($w \in W_0$). (Here $\mathcal{H}_0$ is the subalgebra of $\mathcal{H}$ generated by $T_s$, $s \in S_0$.)

**Definition 2.2.** *The image of an element $h \in \mathcal{H}$ in $End(I_t)$ is called its Laplace transform $\hat{h}(t)$. By the natural identification of $I_t$ with $\mathcal{H}_0$ we will always consider $\hat{h}(t)$ as an element of $End(\mathcal{H}_0)$.*

We list some basic properties of $I_t$:

**Proposition 2.3.**   (a) *The center $\mathcal{Z} \simeq \mathcal{A}^{W_0}$ acts by scalars on $I_t$, according to the character $t$ of $\mathcal{A}$.*
   (b) *Every irreducible module over $\mathcal{H}$ is a quotient of some $I_t$.*
   (c) *If $\hat{h}(t)(T_e) = 0$ for all $t$ in some Zariski-dense subset of $T$, then $h = 0$.*

*Proof.* (a) and (b) are trivial, and (c) follows directly from Theorem 1.9. □

To study the modules $I_t$ more seriously, we need to involve intertwining operators.

### 2.2. The intertwining operators

The action of $W_0$ on $\mathcal{A}$ can "almost" be realized by inner automorphisms of $\mathcal{H}$. Here the word "almost" means that one needs to pass to some formal completion of $\mathcal{H}$ in order to make certain elements of $\mathcal{H}$ invertible. The elements which realize this $W_0$-action by means of conjugation (in a completion of $\mathcal{H}$) are called *intertwining elements*. We do not want to study the formal completions of $\mathcal{H}$ here, since the non-invertibility of the intertwiners in $\mathcal{H}$ is at the heart of the study of representations of $\mathcal{H}$.

Let us first discuss the construction of the intertwining elements themselves. Their definition is closely related to Lusztig's relation 1.10.

**Definition 2.4.** *When $s = s_\alpha \in S_0$ (with $\alpha \in F_0$), we define the intertwining element $R_s$ as follows:*

$$R_s =$$
$$= \begin{cases} (1 - \theta_{-\alpha})T_s + (1 - q_{\alpha^\vee}) & (2\alpha \notin R_{\mathrm{nr}}) \\ (1 - \theta_{-2\alpha})T_s + ((1 - q_{\alpha^\vee/2}q_{\alpha^\vee}) + q_{\alpha^\vee/2}^{1/2}(1 - q_{\alpha^\vee})\theta_{-\alpha}) & (2\alpha \in R_{\mathrm{nr}}) \end{cases}$$
$$= \begin{cases} T_s(1 - \theta_\alpha) + (q_{\alpha^\vee} - 1)\theta_\alpha & (2\alpha \notin R_{\mathrm{nr}}) \\ T_s(1 - \theta_{2\alpha}) + ((q_{\alpha^\vee/2}q_{\alpha^\vee} - 1)\theta_{2\alpha} + q_{\alpha^\vee/2}^{1/2}(q_{\alpha^\vee} - 1)\theta_\alpha) & (2\alpha \in R_{\mathrm{nr}}) \end{cases}$$

**Definition 2.5.** *We introduce the following notation. Let $\pm W_0$ be the group $W_0 \cup -W_0$. Assume $\mathcal{R}$ is a ring equipped with $\pm W_0$-action $\pm w : f \to f^{\pm w}$. Suppose we are given an equivariant set of elements $f_\beta$ of $\mathcal{R}$ (i.e. $f_{\pm w\beta} = f_\beta^{\pm w}$), indexed by the roots $\beta \in R_1$. When $w \in W_0$, we put*

$$f_w := \prod_{\beta \in R_{1,+} \cap w^{-1} R_{1,-}} f_\beta.$$

*We simply write $f$ instead of $f_{w_0}$, where $w_0$ denotes the longest element of $W_0$.*

**Example 2.6.** *When we define $\tilde{q}_\beta = q_{\beta^\vee} q_{2\beta^\vee} \in \mathbf{R}_+$ (with the trivial action of $\pm W_0$), then we get the usual $q(w) = \tilde{q}_w$. We apply the notation of Definition 2.5 frequently to the following cases:*

  (1) $\Delta_\beta := 1 - \theta_{-\beta} \in \mathcal{A}$,
  (2) $c_\beta := c_1(\beta, \cdot) \in \mathbf{C}(T)$,
  (3) $n_\beta := \tilde{q}_\beta \Delta_{-\beta}(\cdot) c_{-\beta}(\cdot) \in \mathbf{C}[T]$ *(hence we may and will consider $n_\beta$ as an element of $\mathcal{A}$, by sending the character $x$ of $T$ to $\theta_x \in \mathcal{A}$)*,
  (4) $D_\beta = n_\beta n_{-\beta}$.

**Lemma 2.7.** *In the situation of Definition 2.5 we have the following simple rules:*

  (1) *When $l(uv) = l(u) + l(v)$, then $f_{uv} = f_u^{v^{-1}} f_v$.*
  (2) *In particular, $f_{w_0 w}^w f_{w^{-1}} = f_{w_0} = f$.*
  (3) $f_{w_0 w}^w f_{w^{-1}}^{-e} = f_{w_0}^w = f^w$.
  (4) $f_w^w = f_{w^{-1}}^{-e}$.

**Theorem 2.8.** *Let $s = s_\alpha \in S_0$ with $\alpha \in F_0$, and let $\beta \in R_1$ be the unique element of $R_1$ which is a positive multiple of $\alpha$. (For convenience, we shall write $s = s_\beta$ with $\beta \in F_1$ in the sequel). Then:*

(1) $R_s \theta_x = \theta_{s(x)} R_s$.
(2) The $R_s$ ($s \in S_0$) satisfy the braid relations corresponding to the Coxeter graph $S_0$.
(3) $R_s^2 \in \mathcal{A}$, and in fact $R_s^2 = D_s$.

*Proof.* (1) is a direct reformulation of Theorem 1.10. (3) is proved by (tedious, but straightforward) direct computation, and is left to the reader. To prove (2) we need to show that if $w = s_1 \ldots s_m$ is a *reduced* word in $W_0$, then the element $R_w := R_{s_1} \ldots R_{s_m} \in \mathcal{H}$ indeed only depends on $w$, and not on the chosen reduced word representing $w$. Let us suppose $R'_w$ is obtained in the same way as $R_w$, but using a different reduced word for $w$. Let $t \in T$ and define $r_w(t) = \hat{R}_w(t)(T_e) \in \mathcal{H}_0$, and similarly $r'_w(t)$. By Proposition 2.3(c), in order to show $R_w = R'_w$, it suffices to show that $r_w(t) = r'_w(t)$ for generic $t \in T$. Now observe that both $r_w(t)$ and $r'_w(t)$ are of the form

(2.1) $$\Delta_w(t^{-1}) T_w + \sum_{u<w} c_{u,x} T_u.$$

On the other hand, both $r_w(t)$ and $r'_w(t)$ are in the $\mathcal{A}$ weight space of $I_t$ with eigenvalue $w(t)$, by (1). Since the dimension of $I_t$ is $|W_0|$, it is clear that these weight spaces all have dimension equal to 1 when $t \in T^{\mathrm{reg}}$. Therefore $r_w(t)$ and $r'_w(t)$ are scalar multiples of each other. But then they have to be equal by 2.1. □

**Definition 2.9.** Let $w \in W_0$, and let $w = s_1 \ldots s_m$ be a reduced word. Define $R_w = R_{s_1} \ldots R_{s_m} \in \mathcal{H}$, and $r_w(t) = \hat{R}_w(t)(T_e) \in \mathcal{H}_0$. Notice that (see Example 2.6):

$$r_{w^{-1}}(wt) r_w(t) = n_w(t) n_w(t^{-1})$$
$$= D_w(t).$$

**Corollary 2.10.** *From the above proof we conclude that for every $w \in W_0$, the weight space $I_t^{wt}$ of weight $wt$ in $I_t$ has dimension 1 when $t$ is regular. The elements $r_w(t) \in \mathcal{H}_0$ (with $w \in W_0$) form a basis of $\mathcal{H}_0$ in this case.*

**Corollary 2.11.** *We have the relation*
$$R_w \theta_x R_{w^{-1}} = D_{w^{-1}} \theta_{wx}.$$

**Definition 2.12.** *When $n_w(t) \neq 0$, we introduce normalized intertwining elements $r_w^0(t) \in \mathcal{H}_0$ by*

(2.2) $$r_w^0(t) = (n_w(t))^{-1} r_w(t).$$

**Corollary 2.13.** *The elements $r_w(w't) \in \mathcal{H}_0$ ($w, w' \in W_0$) are well defined if $n(t)n(t^{-1}) \neq 0$. They satisfy the $W_0$-cocycle relation ($u, v \in W_0$ arbitrary):*

$$(2.3) \qquad r_u^0(vt)r_v^0(t) = r_{uv}^0(t).$$

**Definition 2.14.** *Define, for $w \in W_0$, $t \in T$, and $x \in \mathcal{H}_0 \simeq I_t$,*

$$R(w,t): \quad I_t \quad \to \quad I_{wt}$$
$$x \quad \to \quad xr_{w^{-1}}(wt)$$

**Corollary 2.15.** *This defines intertwining maps such that*
1. *$R(w^{-1}, wt)R(w,t) = D_w(t)\mathrm{Id}_{\mathcal{H}_0}$. In particular, $R(w,t)$ is an isomorphism when $D_w(t) \neq 0$.*
2. *$R(u, vt)R(v, t) = R(uv, t)$ when $l(uv) = l(u) + l(v)$. (Equivalently, $r_u(vt)r_v(t) = r_{uv}(t)$ in this situation.)*
3. *By Corollary 2.13,*

$$r_u(vt)r_v(t) = \frac{n_u(vt)n_v(t)}{n_{uv}(t)} r_{uv}(t).$$

*Therefore we have:*

$$R(u, vt)R(v, t) = \frac{n_u(vt)n_v(t)}{n_{uv}(t)} R(uv, t).$$

**Corollary 2.16.** *When $t$ is regular and $D_{w_0}(t) \neq 0$, then $I_t$ is irreducible. (In fact the following stronger result due to S. Kato [2] holds: $\forall t \in T$, $I_t$ is irreducible if and only if $D_{w_0}(t) \neq 0$. This is much deeper, but we shall not need this fact.)*

2.3. **A nondegenerate sesquilinear pairing**

We identify the module $I_t$ with $\mathcal{H}_0$ by $T_w \otimes 1 \to T_w$ ($w \in W_0$), as we did before. We define a nondegenerate sesquilinear pairing between $I_t$ and $I_{\bar{t}^{-1}}$ by

$$(2.4) \qquad \begin{array}{ccc} I_{\bar{t}^{-1}} \times & I_t & \to \quad \mathbf{C} \\ x & y & \to \quad (x, y) := \tau(x^*y) \end{array}$$

where $x, y \in \mathcal{H}_0$.

**Corollary 2.17.** *This pairing satisfies $(x, yz) = (xz^*, y) = (y^*x, z)$ when $x$, $y$ and $z$ are in $\mathcal{H}_0$.*

**Lemma 2.18.** *The element $r_w(t) \in \mathcal{H}_0$ satisfies $r_w(t)^* = r_{w^{-1}}(w\bar{t}^{-1})$. Hence the intertwining operator $R(w,t)$ satisfies: $R(w,t)^* = R(w^{-1}, w\bar{t}^{-1})$.*

*Proof.* It is sufficient to prove this for $w \in S_0$, and in this case it is an easy direct computation. Use Corollary 2.17. □

**Theorem 2.19.** *The pairing between $I_{\bar{t}^{-1}}$ and $I_t$ respects the $*$ operator: $\forall x \in \mathcal{H}$, $y$, $z \in \mathcal{H}_0$ we have*

$$(\hat{x}^*(\bar{t}^{-1})y, z) = (y, \hat{x}(t)z).$$

*Proof.* When $x \in \mathcal{H}_0$. this is just Corollary 2.17. Hence it is enough to verify the relation when $x \in \mathcal{A}$. It suffices to do so for $t$ in some Zariski-dense set, and thus we may and will assume that $t$ is regular. By Corollary 2.10, $(r_w(t))_{w \in W_0}$ is a basis of eigenvectors of $\hat{\theta}_x(t)$, with eigenvalues $x(wt)$. Likewise, using Proposition 1.12, we have that $T_{w_0} r_{w_0 v}(\bar{t}^{-1})$ is $\hat{\theta}_x^*(\bar{t}^{-1})$ eigenvector with eigenvalue $v\bar{t}(x)$. Hence we just need to show that for $v \neq w$:

$$(T_{w_0} r_{w_0 v}(\bar{t}^{-1}), r_w(t)) = 0.$$

But, using Corollary 2.15 and Lemma 2.18, we have:

$$\begin{aligned}
(T_{w_0} r_{w_0 v}(\bar{t}^{-1}), r_w(t)) &= \tau(r_{v^{-1} w_0}(w_0 v t) T_{w_0} r_w(t)) \\
&= \tau(r_w(t) r_{v^{-1} w_0}(w_0 v t) T_{w_0}) \\
&= const. \tau(r_{w v^{-1} w_0}(w_0 v t) T_{w_0}).
\end{aligned} \quad (2.5)$$

By equation 2.1 this last expression is 0 when $v \neq w$, as was claimed. □

**Proposition 2.20.** *In addition to the proof of the previous theorem we see that, since $l(w) + l(w^{-1} w_0) = l(w_0)$,*

$$(T_{w_0} r_{w_0 v}(\bar{t}^{-1}), r_w(t)) = \delta_{v,w} q(w_0) \Delta(wt). \quad (2.6)$$

### 2.4. Matrix elements of the minimal principal series

We denote by $\mathcal{H}^*$ the linear dual of $\mathcal{H}$, the space of all complex linear functionals on $\mathcal{H}$. Given an endomorphism $\psi \in \text{End}(\mathcal{H}_0)$ we define the corresponding matrix element $E(\psi, t) \in \mathcal{H}^*$ of the minimal principal series module $I_t$ by

$$E(\psi, t)(h) := \text{tr}_{\mathcal{H}_0}(\psi \hat{h}(t)). \quad (2.7)$$

The *character* $\chi(I_t)$ of $I_t$ is by definition the matrix element $E(\text{Id}, t)$. The intertwiners and the pairing discussed above give rise to two functional equations for the matrix elements:

**Proposition 2.21.** (1) $\overline{E(\psi^*, \bar{t}^{-1})(h^*)} = E(\psi, t)(h)$.
(2) $E(R(w, t) \psi R(w^{-1}, wt), wt) = D_w(t) E(\psi, t)$.

With the pairing of the previous subsection at our disposal, we define an important collection of matrix elements of the minimal principal series modules.

**Definition 2.22.** *For $u, v \in W_0$ we define the linear functional $E_{u,v}$ on $\mathcal{H}$ by ($h \in \mathcal{H}$):*

$$(2.8) \qquad E_t^{u,v}(h) = (T_{w_0} r_{w_0 u}(\bar{t}^{-1}), \hat{h}(t) r_v(t)).$$

*In other words, $E_{u,v} = E(\psi_{u,v}(t), t)$ where $\psi_{u,v}(t) \in \mathrm{End}(\mathcal{H}_0)$ is defined by $\psi_{u,v}(t)(h) = (T_{w_0} r_{w_0 u}(\bar{t}^{-1}), h) r_v(t)$.*

**Proposition 2.23.**  (1) *The function $t \to E_t^{u,v}(h)$ is a regular function on $T$, for all $u, v$ and $h$.*
  (2) $E_t^{u,v}(\theta_{x_1} h \theta_{x_2}) = u(t)(x_1) v(t)(x_2) E_t^{u,v}(h)$.
  (3) *When $t$ is regular and $c(t)c(t^{-1}) \neq 0$, $\{E_t^{u,v}\} \subset \mathcal{H}^*$ is a collection of linearly independent linear functionals on $\mathcal{H}$.*
  (4) *If the stabilizer $W_t$ of $t$ is generated by reflections, then the dimension of the linear space $\mathcal{H}_t^* = \{\phi \in \mathcal{H}^* \mid \phi(zh) = z(t)\phi(h), \forall z \in \mathcal{Z}, h \in \mathcal{H}\}$ is equal to $|W_0|^2$. In particular, $\{E_t^{u,v}\} \in \mathcal{H}^*$ is a basis for this space when $t$ is regular and $c(t)c(t^{-1}) \neq 0$.*

*Proof.* (1) is obvious, (2) is straightforward, and (3) follows from (2), if we know in addition that $E_t^{u,v} \neq 0$ when $t$ is regular and $c(t)c(t^{-1}) \neq 0$. This follows from the multiplication rules in Corollary 2.13 for the normalized intertwining elements $r_u^0(t)$ ($u \in W_0$), and equation 2.6.

(4). Note that $\mathcal{H}_t^*$ can be identified with the dual of $\mathcal{H}/(\mathcal{I}_t \mathcal{H})$, where $\mathcal{I}_t$ denotes the ideal in $\mathcal{Z}$ of elements vanishing in $t \in T$. But $\mathcal{H}/(\mathcal{I}_t \mathcal{H}) \simeq \mathcal{H}_0 \otimes \mathcal{A}/(\mathcal{I}_t \mathcal{A})$, thus it is sufficient to show that $\mathcal{A}/(\mathcal{I}_t \mathcal{A})$ has dimension $|W_0|$ if the stabilizer of $t$ is generated by reflections. For $t' \in W_0 t$, let $m_{t'}$ denote the maximal ideal of $t'$ in $\mathcal{A}$, and let $M = \prod_{t' \in W_0 t} m_{t'}$. Note that $\mathcal{I}_t \mathcal{A} \supset M^k$ for suitably large $k$, because the ring $\mathcal{A}/\mathcal{I}_t \mathcal{A}$ is finite dimensional over $\mathbf{C}$. Hence $\mathcal{A}/\mathcal{I}_t \mathcal{A} \simeq \hat{\mathcal{A}}_M/\mathcal{I}_t \hat{\mathcal{A}}_M$, where $\hat{\mathcal{A}}_M$ denotes the $M$-adic completion of $\mathcal{A}$. By the chinese remainder theorem, $\hat{\mathcal{A}}_M \simeq \oplus_{t' \in W_0 t} \hat{\mathcal{A}}_{m_{t'}}$. Furthermore, we have $\mathcal{I}_t \hat{\mathcal{A}}_{m_{t'}} \simeq \mathcal{J}_{t'} \hat{\mathcal{A}}_{m_{t'}}$, where $\mathcal{J}_{t'}$ denotes the ring of $W_{t'}$-invariant elements of $\mathcal{A}$ vanishing in $t'$. Given $f \in \mathcal{J}_{t'}$ and any $k \in \mathbf{N}$, it is possible to find a $F \in \mathcal{I}_t$ such that $ef \in F + m_{t'}^k$ with $e \in \hat{\mathcal{A}}_{m_{t'}}$ invertible. This implies the claim, since $m_{t'}^k$ is contained in both $\mathcal{I}_t \hat{\mathcal{A}}_{m_{t'}}$ and $\mathcal{J}_{t'} \hat{\mathcal{A}}_{m_{t'}}$ if $k$ is sufficiently large.

Since $W_{t'}$ is generated by reflections, we see by the theorem of Chevalley that the dimension of $\hat{\mathcal{A}}_{m_{t'}}/\mathcal{J}_{t'} \hat{\mathcal{A}}_{m_{t'}}$ equals $|W_{t'}| = |W_t|$ for each $t' \in W_0 t$. Whence the result. $\square$

We now list some straightforward consequences of Proposition 2.20, Proposition 2.21 and Corollary 2.15.

**Proposition 2.24.**   (1) $E_t^{u,v}(T_e) = \delta_{u,v} q(w_0) \Delta(ut)$.
(2) $\chi(I_t) = q(w_0)^{-1} \sum_{w \in W_0} \Delta(wt)^{-1} E_t^{w,w}$.
(3) $E_t^{u,v}(h) = \overline{E_{\bar{t}^{-1}}^{w_0 v, w_0 u}(T_{w_0}^{-1} h^* T_{w_0})}$.
(4) $E_t^{u,v} = \frac{n_v(t) n_{uw^{-1}}(wt)}{n_u(t) n_{vw^{-1}}(wt)} E_{wt}^{uw^{-1}, vw^{-1}}$.

**Definition 2.25.** *The element $E_t^{e,e} \in \mathcal{H}^*$ plays a predominant role, and will be denoted by $E_t$.*

**Corollary 2.26.** *The character $\chi(I_t)$ equals $q(w_0)^{-1} \sum_{w \in W_0} \Delta(wt)^{-1} E_{wt}$.*

## 2.5. Macdonald's spherical function

Let $T_0^+$ denotes the central idempotent of $\mathcal{H}_0$ corresponding to the trivial representation $T_w \to q(w)$. In other words,
$$T_0^+ = P_0(q)^{-1} \sum_{w \in W_0} T_w,$$
where $P_0$ is the so called Poincaré polynomial
$$P_0(q) = \sum_{w \in W_0} q(w)$$
of the Weyl-group $W_0$. Macdonald's spherical function $\phi_t \in \mathcal{H}_t^*$ is the matrix coefficient $E(\psi^+, t)$ for the endomorphism
$$\psi^+(h) = P_0(q)(T_0^+, h) T_0^+.$$
Note that the chosen normalization is such that $\phi_t(T_e) = \mathrm{Tr}_{\mathcal{H}_0}(\psi^+) = 1$. It is not difficult to express $\phi_t$ in terms of the basis $E_t^{u,v}$ we introduced in the previous subsection.

**Lemma 2.27.**   (1) $\forall w \in W_0,\ T_0^+ r_w^0(t) = T_0^+$.
(2) $T_0^+ = \frac{q(w_0)}{P_0(q)} \sum_{w \in W_0} c(wt) r_w^0(t)$.

*Proof.*   (1) is a simple direct computation.
(2) Write
$$(2.9) \qquad T_0^+ = \sum_{w \in W_0} b_w(t) r_w^0(t).$$

Using (1) we find that $b_w(t) = b_{w_0}(w_0 w t)$. We compute $b_{w_0}$ by the remark that the only summand on the right hand side of 2.9 that contributes to the coefficient of $T_{w_0}$, is $b_{w_0}(t)r^0_{w_0}(t)$. Now recall that, modulo the space spanned by $T_w$ with $w \neq w_0$,

$$r^0_{w_0}(t) = \frac{\Delta(t^{-1})}{n(t)} T_{w_0}.$$

Comparing coefficients leads to the required result. □

Clearly, $\phi_t$ is completely determined by its values on the subalgebra $\mathcal{H}^+ = T_0^+ \mathcal{H} T_0^+$. This important $*$-subalgebra of $\mathcal{H}$ is isomorphic to the center of $\mathcal{H}$ by the so-called Satake-isomorphism, via the map

(2.10)
$$\begin{aligned} \mathcal{Z} &\to \mathcal{H}^+ \\ z &\to T_0^+ z \end{aligned}$$

It is also clear that the elements $\theta_x^+ = T_0^+ \theta_x T_0^+$ with $x \in X^+$ form a linear basis of $\mathcal{H}^+$. Let us determine the spherical function by computing its values on this basis. The result is the well known formula of Macdonald for the spherical function.

**Theorem 2.28.** *(Macdonald [4]) When $x \in X^+$,*

$$\phi_t(\theta_x^+) = \frac{q(w_0)}{P_0(q)} \sum_{w \in W_0} c(wt) wt(x).$$

*Proof.* Use the above lemma. □

It is now easy to relate $\phi_t$ to the matrix element $E_t$:

**Proposition 2.29.**

$$E_t(T_0^+ h T_0^+) = \frac{q(w_0) n(t^{-1})}{P_0(q)} \phi_t(h).$$

*Proof.*

$$\begin{aligned} E_t(T_0^+ h T_0^+) &= (T_{w_0} r_{w_0}(\bar{t}^{-1}), T_0^+ \hat{h}(t) T_0^+) \\ &= (T_0^+ T_{w_0} r_{w_0}(\bar{t}^{-1}), \hat{h}(t) T_0^+) \\ &= q(w_0) n(t^{-1}) (T_0^+ r^0_{w_0}(\bar{t}^{-1}), \hat{h}(t) T_0^+) \\ &= q(w_0) n(t^{-1}) (T_0^+, \hat{h}(t) T_0^+) \text{ by Lemma 2.27(1)} \\ &= \frac{q(w_0) n(t^{-1})}{P_0(q)} \phi_t(h). \end{aligned}$$

□

The basis $\theta_x^+$ ($x \in X^+$) of $\mathcal{H}^+$ is orthogonal. We need the following (standard) notation to formulate the result. When $x \in X^+$, let $W_x$ denote the stabilizer of $x$ in $W$. Let $W^x$ be the set of shortest length representatives for the right cosets of $W_x$ in $W_0$. It is well known that $q(uv) = q(u)q(v)$ and that $T_{uv} = T_u T_v$ when $u \in W^x$ and $v \in W_x$.

**Definition 2.30.** Let $P_x(q) = \sum_{w \in W_x} q(w)$, and $P^x(q) = \sum_{w \in W^x} q(w)$. Note that $P_0 = P^x P_x$. Likewise, define $T_x^+ = P_x(q)^{-1} \sum_{w \in W_x} T_w$ and $T^{x,+} = P^x(q)^{-1} \sum_{w \in W^x} T_w$. Then $T_0^+ = T^{x,+} T_x^+$. Finally we write $w_x$ for the longest element in $W_x$, and $w^x \in W^x$ for the element such that $w_0 = w^x w_x$.

**Proposition 2.31.** When $x, y \in X^+$, we have:
$$(\theta_x^+, \theta_y^+) = \frac{\delta_{x,y} q(w^x)}{P_0(q) P^x(q)}.$$

*Proof.* Notice that $\theta_x^+$ is a linear combination of elements of the form $T_{u t_x v}$ with $u \in W_0$ and $v \in W^x$. In fact we have:

(2.11)
$$\begin{aligned}
\theta_x^+ &= \delta(-x)^{1/2} T_0^+ T_{t_x} T_0^+ \\
&= \delta(-x)^{1/2} T^{x,+} T_{t_x} T_0^+ \\
&= \frac{\delta(-x)^{1/2} q(w^x)}{P_0(q) P^x(q)} \sum_{u \in W^x, v \in W_0} T_{u t_x v}
\end{aligned}$$

To verify the last line, notice first that the sum on the right hand side is invariant for multiplication by $T_0^+$ on either side. So we only need to check the coefficient of one, suitably chosen element $T_{u t_x v}$. Choose $T_{w^x t_x w_0}$. With the above expression for $\theta_x^+$ at our disposal, it is easy to compute $(\theta_x^+, \theta_y^+) = \delta(-y)^{1/2} (\theta_x^+, T_{t_y})$. $\square$

## 3. Eisenstein series for the Hecke algebra

Consider the vector space of formal sums $\sum_{w \in W} c_w T_w$. Notice that this vector space carries a natural structure of a left and right $\mathcal{H}$ module, since multiplication (on the left or the right) with finite sums is always well defined. Also the trace $\tau$ has a natural extension to the vector space of formal sums. When $\phi = \sum_{w \in W} c_w T_w$, we shall thus define a linear functional $\phi$ on $\mathcal{H}$ by $\phi(h) = \tau(\phi h)$. In other words, we define $\phi(T_w) = q(w) c_{w^{-1}}$. In this way we will identify the algebraic dual $\mathcal{H}^*$ with the vector space of formal sums. We equip $\mathcal{H}^*$ with the weak topology.

Now let us consider, for $t \in T$, the sum
$$\mathcal{E}_t = \sum_{x \in X} t(-x) \theta_x.$$

If this sum is convergent in $\mathcal{H}^*$, it will clearly satisfy:
$$\mathcal{E}_t(\theta_x h) = \mathcal{E}_t(h\theta_x) = t(x)\mathcal{E}_t(h).$$
By Proposition 2.23 this implies that
$$\mathcal{E}_t = f(t) E_t,$$
provided that the left hand side converges. We want to compute the function $f$, but let us first treat the question of convergence.

### 3.1. Convergence in $\mathcal{H}^*$

**Lemma 3.1.** *Let $u, v \in W_0$, and $x \in X$. Then*
  (1) $\tau(T_u \theta_x T_v) = 0$ *unless $x \in Q_-$.*
  (2) $\forall \epsilon > 1$, $\exists C_\epsilon > 0$ *such that $\forall x \in Q_-$,*

$$|\tau(T_u \theta_x T_v)| \leq C_\epsilon \delta_{\epsilon q}^{1/2}(-x). \tag{3.1}$$

*Proof.* (1) When $x \in X^+$ it is clear that $T_u \theta_x T_v$ will be a linear combination of $T_w$ with $w$ in the double coset $W_0 t_x W_0$. Hence this will not have a constant term, unless $x = 0$. For general $x \in X$, we denote by $C_+(x)$ the convex hull of the $W_0$ orbit of $x$, intersected with the cone $x + Q_+$. We claim that for every $x \in X$:

$$T_u \theta_x T_v \in \sum_{y \in X^+ \cap C_+(x)} H_0 \theta_y H_0. \tag{3.2}$$

Let us prove 3.2. When $x \notin X^+$, we can choose a fundamental root $\alpha \in F_0$ such that $(x, \alpha^\vee) < 0$. Apply Lusztig's formula (Theorem 1.10) with $s = s_{\alpha^\vee}$, and we see that we can express $\theta_x$ as follows:

$$\theta_x \in \sum_{y \in C_+(x),\ y \neq x} H_0 \theta_y H_0.$$

When there still are elements $y \in X \backslash X^+$ in this sum, we repeat this procedure for those $y$. After $k$ steps, we have expressed $\theta_x$ as an element in the sum of double coset spaces $H_0 \theta_y H_0$, with $y$ either in the set $X^+ \cap C_+(x)$, or otherwise in $\{y \in C_+(x) \mid (y - x, \rho^\vee) \geq k\}$. By the finiteness of $C_+(x)$, this last set will be empty if $k$ is sufficiently large.
(2) We use induction on the height $-(x, \rho^\vee)$ for $x \in Q_-$. First we choose $N \in \mathbf{N}$ such that $\forall \alpha \in F_0$:

$$(2q_{\alpha^\vee} q_{\alpha^\vee/2} - 1)(\epsilon q_{\alpha^\vee})^{-N} + \frac{\epsilon q_{\alpha^\vee} q_{\alpha^\vee/2}(q_{\alpha^\vee} - 1) + q_{\alpha^\vee} q_{\alpha^\vee/2} - 1}{\epsilon^2 q_{\alpha^\vee}^2 q_{\alpha^\vee/2} - 1} \leq 1. \tag{3.3}$$

We use the usual convention that $q_{\alpha^\vee/2} = 1$ if $\alpha^\vee \notin 2Y$, in which case condition 3.3 reduces to

$$(3.4) \qquad (2q_{\alpha^\vee} - 1)(\epsilon q_{\alpha^\vee})^{-N} + \frac{q_{\alpha^\vee} - 1}{\epsilon q_{\alpha^\vee} - 1} \leq 1.$$

Write $\rho^\vee = \sum_{\alpha \in F_0} l_\alpha \alpha^\vee$. Choose $M \in \mathbf{N}$ such that
$$M > N|F_0|\max\{l_\alpha\}_{\alpha \in F_0}.$$
Consequently, if $-(x, \rho^\vee) \geq M$ then $\exists \alpha \in F_0$ such that $-(x, \alpha^\vee) \geq N$. In order to start the induction, choose $C_\epsilon > 0$ such that equation 3.1 holds $\forall x \in Q_-$ for which $-(x, \rho^\vee) < M$ (a finite subset of $Q_-$). Let $x \in Q_-$ with $-(x, \rho^\vee) \geq M$, and assume by induction that 3.1 holds $\forall y \in Q_-$ such that $(y, \rho^\vee) > (x, \rho^\vee)$. Choose $\alpha \in F_0$ such that $-(x, \alpha^\vee) \geq N$, which is possible by our choice of $M$, and let $s = s_{\alpha^\vee}$. We write (assuming $\alpha^\vee \in 2Y$, the other case being similar and easier):

(3.5)
$$T_u \theta_x T_v =$$
$$= T_u T_s (T_s^{-1} \theta_x T_s) T_s^{-1} T_v$$
$$= T_u T_s (\theta_x \theta_{-\alpha}^{(x, \alpha^\vee)}) T_s^{-1} T_v +$$
$$+ T_u \left( (q_{\alpha^\vee/2} q_{\alpha^\vee} - 1) + q_{\alpha^\vee/2}^{1/2}(q_{\alpha^\vee} - 1)\theta_{-\alpha} \right) \theta_x \frac{1 - \theta_{-\alpha}^{(x, \alpha^\vee)}}{1 - \theta_{-2\alpha}} T_s^{-1} T_v$$

Now we note that
(3.6)
$$T_u T_s = \begin{cases} T_{us} & \text{if } l(us) = l(u) + 1 \\ (q_{\alpha^\vee} q_{\alpha^\vee/2} - 1) T_u + q_{\alpha^\vee} q_{\alpha^\vee/2} T_{us} & \text{if } l(us) = l(u) - 1 \end{cases}$$

and that, similarly,
(3.7)
$$T_s^{-1} T_v = \begin{cases} T_{sv} & \text{if } l(sv) = l(v) - 1 \\ (q_{\alpha^\vee}^{-1} q_{\alpha^\vee/2}^{-1} - 1) T_v + q_{\alpha^\vee}^{-1} q_{\alpha^\vee/2}^{-1} T_{sv} & \text{if } l(sv) = l(v) + 1 \end{cases}$$

Now 3.5 leads directly to the desired result when we use the simple equations 3.6 and 3.7, together with the induction hypothesis, the fact that $-(x, \alpha^\vee) \geq N$, and the inequality 3.3. $\square$

**Corollary 3.2.** $\mathcal{E}_t$ *is weakly convergent if* $\operatorname{Re}(t) < \delta^{-1/2}$. *Here we use the ordering on the space* $T_{rs} = \operatorname{Hom}(X, \mathbf{R}_+)$ *of real characters of* $X$, *given by* $t_1 < t_2 \Leftrightarrow t_1(\alpha) < t_2(\alpha) \forall \alpha \in R_{0,+}$.

*Proof.* On a basis element $T_u \theta_y$ we have to show that
$$\sum_{x \in X} |t(y - x)| |\tau(T_u \theta_x)|$$

is convergent. By Lemma 3.1 it is enough to check the convergence of
$$\sum_{x \in Q_+} \text{Re}(t)(x)\delta^{1/2}_{\epsilon q}(x).$$
This is clear when we choose $\epsilon$ sufficiently small. $\square$

**Definition 3.3.** *The series $\mathcal{E}_t \in \mathcal{H}^*$ will be called "the Eisenstein series" for $\mathcal{H}$.*

## 3.2. Meromorphic continuation of $\mathcal{E}_t$

The formal series $\mathcal{E}_t$ has meromorphic continuation to $T$, and is in fact a rational function. This simple fact is proved in the next lemma.

**Lemma 3.4.** *Recall the notations of Example 2.6. The functional $D(t)\mathcal{E}_t$ (with $\text{Re}(t) < \delta^{-1/2}$) can be written as:*
$$D(t)\mathcal{E}_t = \sum_{w \in W_0} D^w_{w_0 w} R_w \mathcal{E}^w_t R_{w^{-1}},$$
*where*
$$\mathcal{E}^w_t = \sum_{x \in X^+ : w \in W^x} t(-wx)\theta_x.$$
*In particular, $\forall h \in \mathcal{H}$, $t \to D(t)\mathcal{E}_t(h)$ extends to a regular function on $T$.*

*Proof.* This is clear from the relation given in Corollary 2.11, since
$$D(t)\mathcal{E}_t = D\mathcal{E}_t$$
$$= \sum_{w \in W_0} D^w_{w_0 w} D_{w^{-1}} \sum_{x \in X^+ : w \in W^x} t(-wx)\theta_{wx}$$
$$= \sum_{w \in W_0} D^w_{w_0 w} R_w \left( \sum_{x \in X^+ : w \in W^x} t(-wx)\theta_x \right) R_{w^{-1}}.$$
$\square$

**Theorem 3.5.** *Recall Definition 2.25 of the matrix element $E_t$. We have:*
$$D(t)\mathcal{E}_t = \Delta(t^{-1})E_t.$$

*Proof.* As we already remarked in the beginning of this section, it is clear that $E_t$ and $\mathcal{E}_t$ are proportional. In order to find the ratio of proportionality, we project both of them onto Macdonald's spherical function $\phi_t$. In the case of $E_t$ this was done in Proposition 2.29.

Now let us concentrate on $\mathcal{E}_t$. We write

$$D(t)\mathcal{E}_t(\theta_x^+) = D\mathcal{E}_t(\theta_x^+)$$
$$= \tau(T_0^+ D\mathcal{E}_t T_0^+ \theta_x)$$
(3.8)
$$= \sum_{w \in W_0} \tau(T_0^+ D_{w_0 w}^w R_w \mathcal{E}_t^w R_{w^{-1}} T_0^+ \theta_x).$$

Before we continue, it is convenient to introduce some more notation. Let $f$ and $g$ be elements in $\mathcal{H}^*$. We say that $f$ and $g$ are asymptotically equal if there exists an $N \in \mathbf{N}$ such that $f(T_{wt_x}) = g(T_{wt_x})$, $\forall w \in W_0$ and $\forall x \in X$ that satisfy the condition: $\forall \alpha \in R$, $|(x, \alpha^\vee)| \geq N$. This is an equivalence relation on $\mathcal{H}^*$, which we shall denote by $f \sim g$. Note that $\sim$ is respected by the left and right $\mathcal{H}$ module structure of $\mathcal{H}^*$, and that the equivalence class of $0 \in \mathcal{H}^*$ contains $\mathcal{H} \subset \mathcal{H}^*$.

Given $t \in T$, we introduce a right and a left evaluation map, denoted $h \to h(t)$ and $h \to (t)h$ respectively, from $\mathcal{H}$ to $\mathcal{H}_0$. These maps are defined on basis elements by:

$$(T_w \theta_x)(t) = t(x) T_w$$

and

$$(t)(\theta_x T_w) = t(x) T_w.$$

Of course, the right (left) evaluation is the unique homomorphism of left (right) $\mathcal{H}_0$ modules extending the usual evaluation map $\theta_x \to t(x) T_e$ on $\mathcal{A}$.

With these notations it is clear that

$$D_{w_0 w}^w R_w \mathcal{E}_t^w R_{w^{-1}} \sim D_{w_0 w}(w^{-1}t) R_w(w^{-1}t) \mathcal{E}_t^w \left((w^{-1}t) R_{w^{-1}}\right)$$

Moreover, using the Definitions 2.4 and 2.9 we easily find:

$$R_w(w^{-1}t) = r_w(w^{-1}t)$$

and

$$(w^{-1}t) R_{w^{-1}} = r_{w^{-1}}(t^{-1}) \prod_{\beta \in R_{1,+} \cap w R_{1,-}} (-t(\beta)).$$

After this simplification we can deal with the left and right multiplication by $T_0^+$ that occurs in 3.8. We can use the fact that $T_0^+$ is central in $\mathcal{H}_0$, and that, by Lemma 2.7 and Lemma 2.27, we have:

$$T_0^+ r_w(w^{-1}t) = n_{w^{-1}}(t^{-1}) T_0^+$$

and

$$r_{w^{-1}}(t^{-1}) T_0^+ = T_0^+ n_{w^{-1}}(t^{-1}).$$

Now continue equation 3.8, to obtain *for $x \in X^+$ very far from walls*:

$D(t)\mathcal{E}_t(\theta_x^+) =$

$$= \sum_{w \in W_0} \tau(T_0^+ D_{w_0w}^w R_w \mathcal{E}_t^w R_{w^{-1}} T_0^+ \theta_x)$$

$$= \sum_{w \in W_0} \{\prod_{\beta \in R_{1,+} \cap wR_{1,-}} (-t(\beta))\} D_{w_0w}(t) n_{w^{-1}}(t^{-1})^2 \tau(T_0^+ \mathcal{E}_t^w T_0^+ \theta_x)$$

$$= \sum_{w \in W_0} \{\prod_{\beta \in R_{1,+} \cap wR_{1,-}} (-t(\beta))\} \{n_{w_0w}(w^{-1}t) n_{w^{-1}}(t^{-1})\}$$

$$\{n_{w_0w}(w^{-1}t^{-1}) n_{w^{-1}}(t^{-1})\} \sum_{y \in X^+} t(-wy) \tau(\theta_y^+ \theta_x)$$

$$= \frac{q(w_0)}{P_0(q)^2} \sum_{w \in W_0} \{\prod_{\beta \in R_{1,+} \cap wR_{1,-}} (-t(\beta))\} n^w(t) n(t^{-1}) (w_0 w^{-1} t)(x)$$

(using Lemma 2.7 and Proposition 2.31)

$$= \frac{q(w_0)^3}{P_0(q)^2} \Delta(t) \Delta(t^{-1}) c(t) \sum_{w \in W_0} c(wt) wt(x)$$

$$= \frac{q(w_0)^2}{P_0(q)} \Delta(t) \Delta(t^{-1}) c(t) \phi_t(\theta_x^+)$$

$$= \frac{q(w_0)}{P_0(q)} \Delta(t^{-1}) n(t^{-1}) \phi_t(\theta_x^+)$$

Comparing this result with Proposition 2.29 finally gives the desired result. □

**Remark 3.6.** *Surprisingly, the result of the above computation is correct for all strongly dominant $x$ in $D(t)\mathcal{E}_t(\theta_x^+)$. This is better than one would expect, at first sight, from this method of computation!*

### 3.3. Proof of the main theorem

By Theorem 3.5 we may write, when $\text{Re}(t) < \delta^{-1/2}$:

$$(3.9) \quad \mathcal{E}_t = \sum_{x \in X} t(-x) \theta_x = \left(\frac{E_t}{q(w_0)\Delta(t)}\right) \left(\frac{1}{q(w_0)c(t^{-1})c(t)}\right)$$

The main Theorem 1.14 follows by taking the trace $\tau$ in formula 3.9, using Proposition 2.24.

## 3.4. A disintegration of the trace

The formula for the trace we have now derived is only a simple formal consequence of formula 3.9. Formula 3.9 is itself of fundamental interest, since it is the starting point for the spectral analysis of $\mathcal{H}$. This may be illustrated by the following obvious disintegration formula : if $t_0 \in T_{rs} \in \text{Hom}(X, \mathbf{R}_+)$, the split real form of $T$, with $t_0 < \delta^{-1/2}$, then

$$
(3.10) \quad \begin{aligned} T_e &= \sum_{x \in X} \int_{t \in t_0 T_c} t(-x) \theta_x dt \\ &= \int_{t \in t_0 T_c} \mathcal{E}_t dt. \end{aligned}
$$

Thus, when we use 3.9, and interpret formula 3.10 as a weak integral of a function with values in $\mathcal{H}^*$, we obtain:

**Theorem 3.7.**

$$
(3.11) \quad \tau = \int_{t \in t_0 T_c} \left( \frac{E_t}{q(w_0) \Delta(t)} \right) \frac{dt}{q(w_0) c(t^{-1}) c(t)}.
$$

Here $T_c = \text{Hom}(X, S^1)$, the compact form of the algebraic torus $T$, and $dt$ denotes the holomorphic $n$-form on $T$ which restricts to the normalized Haar measure on $T_c$. This formula might be called "the Laplace inversion formula" for $\mathcal{H}$. In order to refine this formula to the level of a spectral resolution of $\tau$, we need to carry out a contour shift, sending $t_0$ to $e \in T_{rs}$. The resulting formula will be an integral of tempered characters of $\mathcal{H}$, against a positive measure. This refinement of Theorem 3.7 requires several techniques, and will be discussed elsewhere.

Macdonald's spectral resolution [4] of the trace function of the algebra $\mathcal{H}^+$ can be recovered easily from Theorem 3.7, by applying the projection onto the spherical function as in Theorem 3.5. We find:

$$
(3.12) \quad \begin{aligned} \tau|_{\mathcal{H}^+} &= \int_{t \in t_0 T_c} T_0^+ \mathcal{E}_t T_0^+ dt \\ &= \frac{1}{P_0(q)} \int_{t \in t_0 T_c} \phi_t \frac{dt}{c(t^{-1})} \\ &= \frac{1}{P_0(q)} \int_{t \in T_c} \phi_t \frac{dt}{c(t^{-1})} \end{aligned}
$$

We can send $t_0$ to $e$ in $T_{rs}$ without meeting poles of the integrand, because the projection of $\mathcal{E}_t$ onto $\phi_t$ created a "window" in the set of poles. This step is crucial, and explains the simplicity of the spherical harmonic analysis. In particular there is only continuous spectrum in this case.

We resume the computation of 3.12:

$$\tau|_{\mathcal{H}^+} = \frac{1}{|W_0|P_0(q)} \sum_{w \in W_0} \int_{t \in T_c} \phi_t \frac{dt}{c(wt^{-1})}$$

$$= \frac{1}{|W_0|P_0(q)} \int_{t \in T_c} \phi_t \left( \sum_{w \in W_0} c(wt) \right) \frac{dt}{c(t)c(t^{-1})}$$

$$= \frac{1}{|W_0|q(w_0)} \int_{t \in T_c} \phi_t \frac{dt}{c(t)c(t^{-1})}.$$

More difficult is the formula that one obtains after projection onto the anti-spherical function $\phi_t^-$ defined by $\phi_t^-(h) = P_0(q)(T_0^-, \hat{h}(t)T_0^-)$, with $T_0^-$ the central idempotent of $\mathcal{H}_0$ corresponding to the sign representation $T_w \to (-1)^{l(w)}$ of $\mathcal{H}_0$. Here one obtains:

(3.13)
$$\tau|_{\mathcal{H}^-} = \int_{t \in t_0 T_c} T_0^- \mathcal{E}_t T_0^- \, dt$$

$$= \frac{1}{P_0(q^{-1})} \int_{t \in t_0 T_c} \phi_t^- \frac{dt}{c(t)}.$$

This time, the window in the set of poles of the integrand is in the direction of the positive chamber of $T_{rs}$, so we cannot avoid the poles when we shift the contour. Nonetheless, the structure of the set of poles which one has to deal with, is essentially simpler than in the case of Theorem 3.7. The anti-spherical case was dealt with in [1]. It gives insight in the anti-spherical Plancherel measure, in particular in the formal degrees of the discrete spectrum (with anti-spherical vector!).

We expect that it is possible to find similar information for the full tempered spectrum, starting from Theorem 3.7.

# References

[1] G.J. Heckman, E.M. Opdam, *Harmonic analysis for affine Hecke algebras*, Current Developments in Mathematics (S.-T. Yau, editor), 1996, Intern. Press, Boston.

[2] S. Kato, *Irreducibility of principal series representations for Hecke algebras*, Inv. Math. **87**, (1983), 929-943.

[3] G. Lusztig, *Affine Hecke algebras and their graded version*, J. Amer. Math. Soc. **2** (1989), no. 3, pp. 599–635.

[4] I.G. Macdonald, *Sperical functions on a group of p-adic type*, Publ. Ramanujan Institute **2** (1971).

[5] R. Steinberg, *On a theorem of Pittie*, Topologie **14** (1975), 173-177.

Korteweg-De Vries Institute for Mathematics, Universiteit van Amsterdam, Plantage Muidergracht 24, 1018 TV Amsterdam, The Netherlands.

*E-mail address*: `opdam@wins.uva.nl`